\documentclass[12pt]{article}
\usepackage[utf8]{inputenc}
\usepackage{amsfonts,amsmath,amssymb,amsthm,graphicx}
\input{xy}
\xyoption{all}
\usepackage{graphicx}
\newtheorem*{theorem2}{Conway-Gordon-Sachs Theorem}
\newtheorem*{theorem}{Linear Conway-Gordon-Sachs Theorem}
\newtheorem*{th3}{Sachs Theorem}
\newtheorem*{Pro}{Remark 4$'$}
\newtheorem*{lemma3''}{Lemma 3$''$}
\newtheorem*{lemma2'}{Lemma 3$'$}
\newtheorem*{lemma1'}{Lemma 1$'$}
\theoremstyle{mydefinition}
\newtheorem{lemma}{Lemma}
\newtheorem{example}[lemma]{Example}
\newtheorem*{example'}{Example 2$'$}
\newtheorem{remark}[lemma]{Remark}

\def\Int{\mathop{\fam0 Int}}

\def\R{{\mathbb R}} \def\Z{{\mathbb Z}}

\long\def\comment#1\endcomment{}

\begin{document}
\begin{center}
\textbf {Alternative proofs of the Conway-Gordon-Sachs Theorems}
\end{center}
\begin {center}
Arseny Zimin
\end{center}
\begin{abstract}In this paper we present 
new proofs of the Conway-Gordon-Sachs and Sachs Theorems on the linked cycles in graphs embedded in $\R^3$.
We reduce these theorems to certain property of graphs 
mapped to the plane.
\end{abstract}

By a {\it space} we always mean 3-dimensional space $\R^3$. Points in space are in \textit{general position}, if no four of them are in one plane. By a {\it triangle} we always mean a closed broken line which has three distinct vertices.
Two triangles in space whose six vertices are in general position are \textit{linked}, if 
the first triangle intersects the convex hull of the second triangle exactly at one point.

\begin{theorem}
Assume that six points in space are in general position. Then there exist two linked triangles with the vertices at these points.
\end{theorem}




A {\it spatial polygon} is a non-self-intersecting broken line (not necessarily closed) in space. Let $a$ and $b$ be disjoint spatial polygons.
For a point $A$ in space denote by $A*a$ the union $\bigcup\limits_{X\in a}AX$ of segments.
A point $A$ is {\it in general position} to $a$ and $b$, if

$\bullet$ no vertex of $b$ belongs to $A*a$ and

$\bullet$ if $X$ is either a vertex of $a$ or a point of $a$ such that $a\cap\Int AX\ne\emptyset$,
then $b\cap AX=\emptyset$.

{\it Remark.} If a point $A$ is in general position to disjoint spatial polygons $a$ and $b$,
then the set $(A*a)\cap b$ is a finite.

Closed disjoint spatial polygons $a$ and $b$ are called {\it linked modulo 2}, if there is a point $A$ in general position to $a$ and $b$ such that $|(A*a)\cap b|\equiv 1\mod2$.
From now on we will always say that two spatial polygons are {\it linked} instead of {\it `linked modulo 2'}.
\comment
Consider a closed broken line $a$ in space. A finite set of triangles
in space is called a {\it membrane spanned by $a$}, if the following conditions hold:

$\bullet$ Any side of $a$ is a side of exactly one triangle from this set;

$\bullet$ Any side of a triangle from this set that is not a side of $a$ belongs to exactly two triangles from this set.

Two disjoint non-self-intersecting closed broken lines $a$ and $b$ in $3$-dimensional space are \textit{linked}, if there exists a membrane, denote it by $A$, spanned by $a$
such that

$\bullet$ $b$ do not intersect outline of any triangle from $A$ and no vertex of $b$ belongs to some triangle from $A$;

$\bullet$ the number of triangles from $A$ that intersect $b$ is odd.
\endcomment


\begin{figure}
[h]\centering
\includegraphics[width=2.5cm]{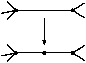}
\caption{subdivision of an edge}
\end{figure}

Operation of {\it subdivision} of an edge of a graph is shown in Fig. 1. Two graphs are called {\it homeomorphic}, if one graph can be obtained from the other one by certain number of subdivisions of edges and inverse operations.
A {\it piecewise-linear embedding} of a graph $G$ into space is a linear embedding of a graph homeomorphic to $G$ into space.

Denote by $K_n$ the complete graph on $n$ vertices.
Denote by $K_{n,n}$ the complete bipartite graph each of whose parts has $n$ vertices.
\begin{theorem2}
Assume that the graph $K_6$ is piecewise-linearly embedded in space. Then there exist two linked cycles of length 3 in this graph.
\end{theorem2}

\begin{th3}
Assume that the graph $K_{4,4}$ is piecewise-linearly embedded in space. Then there exist two linked cycles of length 4 in this graph.
\end{th3}

These theorems are classical results of Ramsey link theory, see more in \cite{PS}. Some applications of Ramsey link theory are presented in, e.g., \cite{SS14}.

We reduce these theorems to results for the plane (see Lemma 3, Lemma 3$'$ and Lemma 3$''$ below).

The original proof of the Conway-Gordon-Sachs Theorem [P06, Th. 1.10]
has two steps. In the first step it is proved that change of the embedding of the graph $K_6$ into space does not change
the parity of number of pairs of linked cycles in this graph.
\footnote{The proof of the first step uses either one of the following two facts:
\newline
$(1)$ any closed spatial polygon and 2-dimensional sphere in space that are in general position intersect at an even number of points;
\newline
$(2)$ any two piecewise-linear embeddings of the graph $K_6$ in space are related by isotopy and passing edges one through another.
\newline
Lemma 1' below follows from the first fact.}
In the second step one constructs an example [P06, Fig. 13] when this number is odd.

In our proof we have planar analogues of these steps (see proofs of Lemma 3 and Lemma 3$'$, Example 2 and Example 2$'$ below).
These analogues are easier, in particular, construction of planar example rather than spatial example is much easier.

The main idea of our proof is shown in the next section. The more technical proof of the Conway-Gordon-Sachs Theorem is given in the section `Proof of the Conway-Gordon-Sachs Theorem'. Analogous proof of the Sachs Theorem is given in the end of the paper.

There is a shorter unpublished proof of the {\it linear Conway-Gordon-Sachs Theorem} invented by Alexander Shapovalov.
That proof does not generalise to the proof of the {\it piecewise-linear} case.

\comment
\begin{figure}
\centerline{\includegraphics[width=2cm]{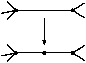}}
\caption{subdivision of an edge}
\end{figure}
\endcomment

\comment

Denote by $K_{n,n}$ the complete bipartite graph each of whose parts has $n$ vertices.
\begin{th3}
Assume that the graph $K_{4,4}$ is piecewise-linearly embedded in $3$-dimensional space. Then there exist two linked cycles of length 4 in this graph.
\end{th3}
\endcomment

\bigskip
{\bf Proof of the linear Conway-Gordon-Sachs Theorem.}

In the proof we use the following known lemmas whose proofs are presented for completeness.

\begin{figure}
\centerline{\includegraphics[width=7cm]{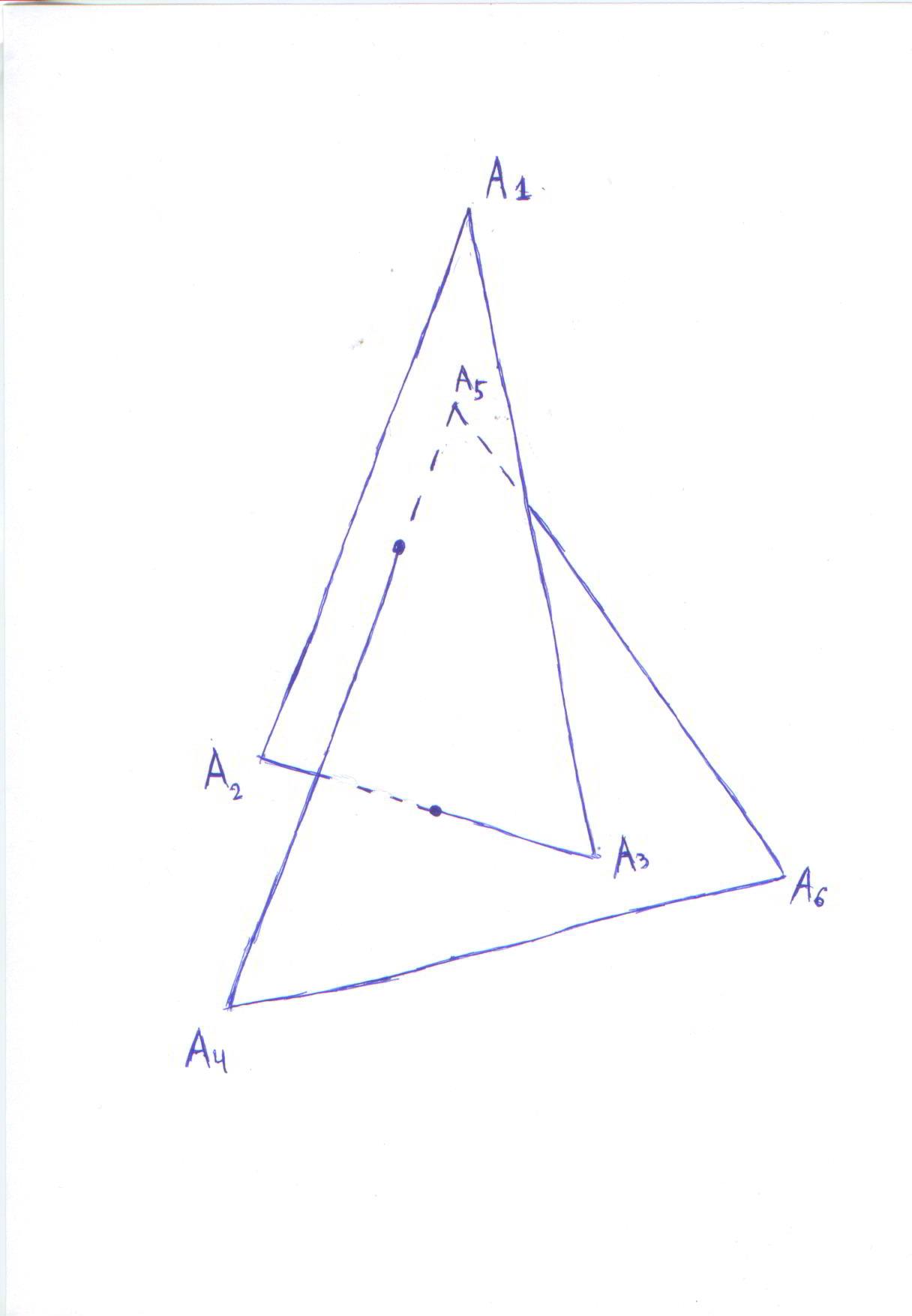}}
\caption{To Lemma \ref{sphere}}
\end{figure}

Let $a,b$ be segments in space. The segment $a$ is \textit{higher} than the segment $b$ looking from
point $O$, if there exists a half-line $OX$ with the endpoint
$O$ that intersects the segment $a$ at a point $A:=a\cap OX$, the segment $b$ at a point $B:=b\cap OX$, $A\not=B$, and $A$ is in the segment
$OB$.

\begin{lemma}[see Figure 2]\label{sphere}
Assume that vertices of triangles $A_1A_2A_3$ and $A_4A_5A_6$
are in general position in space.
If there is a unique side of $A_4A_5A_6$
that is higher than $A_2A_3$ looking from $A_1$,
then these triangles are linked.
\end{lemma}







\smallskip

\textit{Proof.} From the definition of the notion {\it higher} and since the points $A_1,...,A_6$ are in general position, it follows that the number of those sides of the triangle $A_4A_5A_6$ that are higher than $A_2A_3$ is equal to the number of intersection points of the triangle $A_4A_5A_6$ with the convex hull of the triangle $A_1A_2A_3$.
Then the hypothesis implies that the triangle $A_4A_5A_6$ intersects the convex hull of the triangle $A_1A_2A_3$ at a unique point.
{\it QED}

\comment
For any three points $X,Y,Z$ in space denote by $conv(X,Y,Z)$ the convex hull of these point, by $XYZ$ the triangle with vertices $X,Y,Z$ From the definition of the notion {\it higher}, it follows that the number of those sides of the triangle $A_4A_5A_6$ that are higher than $A_2A_3$ is equal to $|A_4A_5A_6\cap (conv(A_1A_2A_3)-A_1A_2A_3)|$. 
Since points $A_1,...,A_6$ are in general position, $A_4A_5A_6\cap conv(A_1,A_2,A_3)=A_4A_5A_6\cap (conv(A_1A_2A_3)-A_1A_2A_3)$.
Then the hypothesis implies that the triangle $A_4A_5A_6$ intersects the convex hull of points $A_1,A_2,A_3$ at a unique point.
\endcomment

\smallskip




\comment
\begin{lemma}\label{vankam}[Sk]
Let a collection $f$ of five points in general position in the plane be given. Then the sum modulo 2 of numbers of intersection points of the segments $AB$ and $CD$ for all unordered pairs
$\{\{A,B\},\{C,D\}\}$ of disjoint two-element subsets $\{A,B\},\{C,D\}\subset f$:

$$v(f):=\sum\{|AB\cap CD|\ :\ \{\{A,B\},\{C,D\}\}\subset{f\choose 2},\ \{A,B\}\cap\{C,D\}=\emptyset  \}\mod2.$$

is equal to 1.
\end{lemma}
\endcomment

Points in the plane are in {\it general position}, if no three of them lie on one line.

\smallskip

{\it Definition of the van Kampen invariant.}
Let $f$ be a set of five general position points in the plane.
For any four distinct points $A,B,C,D\in f$ the segments $AB$ and $CD$ either are disjoint or have a unique common point.
Let {\it the van Kampen invariant} $v(f)\in\Z_2$ be the sum modulo 2 of numbers of intersection points of all the unordered pairs of segments with the vertices at the set $f$ and without common vertices:
$$v(f):=\sum\{|AB\cap CD|\ :\ \{\{A,B\},\{C,D\}\}\subset{f\choose 2},\ \{A,B\}\cap\{C,D\}=\emptyset  \}\mod2.$$

\begin{example}\label{pla-exa}
For the collection $f_0$ of the vertices of a regular pentagon $v(f_0)=1$.
\end{example}

\begin{lemma}\label{vankam}

[Sk]
The van Kampen invariant of each collection of five general position points in the plane is 1.
\end{lemma}

For the proof we need the following simple result.

\begin{remark}\label{0-even}
If the 6 vertices of two triangles in the plane are in general position, then these triangles intersect each other at an even number of points.
\end{remark}
{\it Proof (cf. \cite{BE01}, \S1.5).} Let $A,B,C,D,E,F$ be 6 general position points in the plane. Note that the intersection of the triangle $ABC$ with the convex hull of the triangle $DEF$ is the union of a finite  number of
non-closed broken lines.
The set $ABC \cap DEF$ is formed by the endpoints of these broken lines.
A broken line has 2 endpoints.
Thus $|ABC \cap DEF|$ is even.
{\it QED}

\smallskip

{\it Proof of Lemma \ref{vankam}.}
By Example \ref{pla-exa} it suffices to prove that $v(f)$ does not depend on a set $f$ of five general position points in the plane.
For this, it suffices to prove that if we change one point keeping the remaining four fixed, so that new five points are in general position, then the van Kampen invariant is not changed.
Assume that $V\in f$, $V'\not\in f$ and $f':=(f-\{V\})\cup\{V'\}$ is a general position set.

{\it Proof that $v(f)=v(f')$ when $f\cup\{V'\}$ is a general position set.}
For each $A\in f-\{V\}$ denote by $\triangle_A$ the triangle with the vertices from $f-\{A,V\}$.
Then
$$v(f')-v(f)\overset{1}=\sum\limits_{A\in f-\{V\}}(|VA\cap \triangle_A| -|V'A\cap \triangle_A|) \overset{2}=\sum\limits_{A\in f-\{V\}}|VV'\cap \triangle_A|\overset{3}= 0 \mod 2,\quad\text{where}$$
\ \quad
$\bullet$ the first equality is clear;

$\bullet$ the second equality holds because $|VV'A\cap\triangle_A|$ is even for each $A\in f-\{V\}$ by Remark \ref{0-even};

$\bullet$ the third equality holds because for any two distinct points $M,N\in (f-\{V\})$ the segment $MN$ is contained in exactly two triangles with the vertices from $f-\{V\}$, so the number $|VV'\cap MN|$ `appears twice' in the left-hand sum of the third equality.

{\it Proof that $v(f)=v(f')$ in general.}
There exists a point $V''$ such that $f\cup \{V''\}$ and $f'\cup\{V''\}$ are general position sets.
Then $v(f)=v((f-\{V\})\cup\{V''\})=v(f')$ by the previous case. {\it QED}

\bigskip

\textit{Proof of the linear Conway-Gordon-Sachs Theorem.}
Suppose that points $A_1,A_2,A_3,A_4,A_5,A_6$
are in general position in space.
We may assume that $A_1$ is the unique point among our six points whose first coordinate $a$ is maximal.
Denote by $E$ the set of segments joining pairs of points $A_2,A_3,A_4,A_5,A_6$ and by $\widetilde{E}$ the set of ordered pairs $(e,e')$ of disjoint segments from $E$.
For any ordered pair $(e',e)\in \widetilde{E}$
let the number $lk(e',e)\in \Z_2$ be

$$
lk(e',e):=\begin{cases}
1, &\text{if $e'$ is higher than $e$ looking from $A_1$}\\
0, &\text{otherwise}
\end{cases}.
$$
The number $lk(e',e)$ depends on the point $A_1$ but we omit this from the notation.
For any segment
$e\in E$ denote by

$\bullet$ $A_1e$ the triangle with the vertices at $A_1$ and the endpoints of $e$;

$\bullet$ by $\triangle_e$ the triangle whose vertices are three of our six points other than vertices of the triangle $A_1e$;

$\bullet$ by $lk(\triangle_e, e)$ the number modulo 2 of sides of the triangle $\triangle_e$ that are higher than $e$:

$$lk(\triangle_e,e):=\sum\limits_{e'}{lk(e',e)} \mod 2.$$

Here the summation is over sides $e'$ of $\triangle_e$.

\begin{remark}\label{odd}
For any segment $e\in E$, if
$lk(\triangle_e,e)=1$,
then triangles $\triangle_e$ and $A_1e$ are linked.
\end{remark}
{\it Proof.} Assume that $lk(\triangle_e,e)=1$ for some segment $e\in E$.
Then either exactly one side of $\triangle_e$ is higher than $e$ or any side of $\triangle_e$ is higher than $e$. Assume that
any side of $\triangle_e$ is higher than $e$.
Then from the definition of the notion {\it higher}, it follows that any side of $\triangle_e$ intersects the convex hull of the triangle $A_1e$.
Then these triangles are in one plane and our six points are not in general position. Contradiction. Then exactly one side of $\triangle_e$ is higher than $e$ and
the remark follows from Lemma \ref{sphere}. {\it QED}

\smallskip

{\it Continuation of the proof of the linear Conway-Gordon-Sachs Theorem.} Consider the plane $\alpha$
given by the equation $x=b$, where $b$ is slightly smaller than $a$. Points $A_2,...,A_6$ are in one half-space with respect to $\alpha$ and the point $A_1$ is in the other one.
Let $\pi:\R^3-\{A_1\}\to \alpha$ be the central projection with the center $A_1$. We have

$$\sum\limits_{e\in E}{lk(\triangle_e,e)}\overset{1}=
\sum\limits_{(e',e)\in \widetilde{E}}{lk(e',e)}\overset{2}= v(f)\overset{3}= 1 \mod2,\quad\text{where}$$

$\bullet$ the first equality follows from the definition of $lk(\triangle_e,e)$;

$\bullet$ the second equality follows from Remark \ref{equality} below;

$\bullet$ the third equality is Lemma \ref{vankam}.

Hence for some segment $e\in E$ we have $lk(\triangle_e,e)=1$. Then by Remark \ref{odd} the triangles $\triangle_e$ and $A_1e$ are linked. {\it QED}

\begin{remark}\label{equality}
For any two disjoint segments $e,e'\in E$ we have $lk(e,e') + lk(e',e)\equiv|\pi(e)\cap \pi(e')|\mod 2$.
\end{remark}
{\it Proof.} 
Since points $A_1,...,A_6$ are in general position in space, points $\pi (A_2),..., \pi (A_6)$ are in general position in the plane $\alpha$. Then for any two disjoint segments $e,e'\in E$ the segments $\pi(e)$ and $\pi(e')$ either are disjoint or intersect each other at a unique point.
In the first case the remark holds because for any two disjoint segments $e,e'\in E$, if one of the segments $e$ and $e'$ is higher than the other, then $\pi(e)$ intersects $\pi(e')$.
This fact follows from the definition of the notion higher.
In the second case the remark holds because for any two disjoint segments $e,e'\in E$, if $\pi(e)$ intersects $\pi(e')$, then exactly one of the segments $e$ and $e'$ is higher than the other. This fact follows from the definition of the notion higher and because points $A_2,...,A_6$ are in one half-space with respect to $\alpha$ and the point $A_1$ is in the other one. {\it QED}

\bigskip

{\bf Proof of the Conway-Gordon-Sachs Theorem.}

In the proof we use the following definitions and known lemmas. Some of these lemmas are proved here for completeness.





A set of broken lines in the plane is in {\it general position}, if the following conditions hold:

$\bullet$ no three sides of broken lines from this set have a common interior point;


$\bullet$ no vertex of a broken line from this set lies inside a side of a broken line from this set;

$\bullet$ if two sides of these broken lines have a common vertex, then either they are adjacent sides of one broken line from this set or they are end-sides of two different non-closed broken lines from this set and the common vertex of these sides is a common endpoint of these two broken lines.

Piecewise-linear map $f:G\to\R^2$ of a graph $G$ is called a {\it general position map}, if the set of images of edges of $G$ is in general position.






A plane is {\it in general position} to a graph piecewise-linearly embedded in space, if the orthogonal projection of this graph onto this plane is a general position map.


\begin{lemma}\label{gen}
For any finite graph piecewise-linearly embedded in space
there exists a plane in general position to this graph.
\end{lemma}

{\it Proof.}  The dimension of the set of all directions of orthogonal projections onto planes is equal to the dimension of the projective plane, i. e. to 2. The set of `forbidden' directions of orthogonal projections onto planes is
$(\bigcup\limits_{\{A,e\}}{M_{A,e}})\cup(\bigcup\limits_{\{f,g,h\}}{M_{fgh}})$, where

$\bullet$ $A$ is a vertex of a spatial polygon that is an edge of the graph, $e$ is a side of a spatial polygon that is an edge of the graph;

$\bullet$ $M_{A,e}$ is the set of lines $AX$, where $X$ is a point on $e$;

$\bullet$ $f,g,h$ are sides of spatial polygons that are edges of the graph and $M_{fgh}$ is the set of lines that intersect the union of segments $f,g,h$ exactly at 3 points.



For any vertex $A$ and side $e$ the set $M_{A,e}$ is 1-dimensional because any line from this set depends only on point on $e$. For any three sides $f,g,h$ the set $M_{fgh}$ is at most 1-dimensional because

$\bullet$ if all three sides $f,g,h$ are in one plane, then this set is a subset of projective line;

$\bullet$ if some two of sides $f,g,h$ are in one plane and the third side is not in this plane, then there is at most one line in this set;

$\bullet$ if no two of sides $f,g,h$ are in one plane, then through any point on $f$ passes at most one line from $M_{fgh}$.

Since the graph is finite, $\bigcup\limits_{\{A,e\}}{M_{A,e}}$ and $\bigcup\limits_{\{f,g,h\}}{M_{fgh}}$ are unions of finite numbers of sets. Then the lemma holds. {\it QED}
\comment
Let $\pi$ be a plane in space and $a,b$ segments in space. The segment $a$ is lower than the segment $b$, if there exists a point $O\in\pi$ and a half-line $OX$ with the endpoint
$O$ that is orthogonal to $\pi$ and intersects the segment $a$ at a point $A:=a\cap OX$, and the segment $b$ at a point $B:=b\cap OX$, so that $A\not=B$ and $A$ is in the segment
$OB$.
Define what it means for the first segment to be \textit{higher} than the second segment analogously to the above definition of notion {\it higher} but replacing 'central projection' with 'orthogonal projection'.
\endcomment


\smallskip

Assume that we have orthogonal projection of disjoint spatial polygons $e',e$ onto a plane that is in general position to these polygons and it is shown which of the polygons passes higher than the other at the intersection points of the projections. Let
$$lk(e',e)\in \Z_2$$
be the number modulo 2 of intersection points on the projection at which the polygon $e'$ passes higher than the polygon $e$.

The number $lk(e',e)$ depends on the plane but we omit this from the notation.

\begin{lemma1'}\label{link}
Assume that we have orthogonal projection of two closed disjoint spatial polygons $e',e$ onto a plane that is in general position to these polygons.
If $lk(e',e)=1$, then these polygons are linked.
\end{lemma1'}
This lemma is a generalization of Lemma \ref{sphere}. We do not prove it here.

\smallskip

For any graph $G$ denote by $V(G)$ the set of vertices of $G$, by $E(G)$ the set of edges of $G$ and by $\widetilde{E}(G)$ the set of ordered pairs of disjoint edges of $G$.

\smallskip

{\it Definition of the van Kampen invariant.}
Let $G$ be a graph and $f:G\to\R^2$ a general position map. For any two edges $e,e'$ of $G$ the broken lines $f(e)$ and $f(e')$ intersect at a finite number of points.
Let {\it the van Kampen invariant} $v(f)\in\Z_2$ be the sum modulo 2 of numbers of intersection points of the broken lines $f(e)$ and $f(e')$ for all unordered pairs $\{e,e'\}$ of disjoint edges of $G$:
$$v(f):=\sum\{|f(e)\cap f(e')|\ :\ \{e,e'\}\subset E(G),\ e\cap e'=\emptyset  \}\mod2.$$

\smallskip


\begin{example'}\label{pla-exa'}
For the map $f_0:K_5\to\R^2$ such that $f_0(K_5)$ is a regular pentagon with diagonals $v(f_0)=1$.
\end{example'}

\begin{lemma2'}\label{vankam2}
[Sk] For any general position map $f:K_5\to\R^2$ we have $v(f)=1$.



\end{lemma2'}

This lemma is a generalization of Lemma \ref{vankam}. See the proof below.

\smallskip
\textit{Proof of the Conway-Gordon-Sachs Theorem.}
Assume that the graph $K_6$ is piecewise-linearly embedded in space.
By Lemma \ref{gen} there exists a plane $\alpha$ in general position to $K_6$.
Denote by $A$ one of the vertices of the graph $K_6$. Denote $K_5:=K_6-A$. Let $f:K_5 \to \alpha$ be the orthogonal projection.
For any edge $e$ of $K_5$ denote by 

$\bullet$ $E_1,E_2$ the vertices joined by $e$;

$\bullet$ by $Ae$ the cycle in $K_6$ on the vertices $A,E_1,E_2$;

$\bullet$ by $\triangle_{e}$ the cycle in $K_5$ on three vertices other than $E_1$ and $E_2$.

From now on in any sum, if the limits of the summation are not written, the sum is over edges $e$ of the graph $K_5$.
We have

$$\sum{lk(\triangle_{e}, Ae)}\overset{1}=\sum{lk(\triangle_{e},AE_1)}+\sum{lk(\triangle_{e}, AE_2)}+ \sum{lk(\triangle_{e}, e)}\overset{2}=$$
$$=\sum{lk(\triangle_{e}, e)}\overset{3}=\sum\limits_{(e',e)\in \tilde{E}(K_5)}{lk(e',e)}\overset{4}= v(f)\overset{5}= 1 \mod 2.$$

\comment
$$\sum{lk(abc,\triangle_{bc})}\overset{1}\equiv\sum{lk(ab,\triangle_{bc})}+\sum{lk(ac,\triangle_{bc})}+ \sum{lk(bc,\triangle_{bc})}\overset{2}\equiv$$
$$\equiv\sum{lk(bc,\triangle_{bc})}\overset{3}\equiv\sum{\sum\limits_{e\in E(\triangle_{bc})}{lk(bc,e)}}\overset{4}\equiv\sum\limits_{(e,e')\in \widetilde{E}(K_5)}{lk(e,e')}\overset{5}\equiv v(f)\overset{6}\equiv 1 \mod 2$$
\endcomment




$\bullet$ First and third equalities are clear.

$\bullet$ Let us prove the second equality. For any edge $e$ of $K_5$ we have $lk(\triangle_{e}, AE_1)=\sum\limits_{e'\in E(\triangle_{e})}{lk(e', AE_1)}$.
Any edge $e'$ of the graph $K_5-E_1$ is contained in exactly two cycles of length 3 in $K_5-E_1$. Then the number $lk(e', AE_1)$ `appears twice' in the sum $\sum{lk(\triangle_{e}, AE_1)}$. Therefore this sum is equal to 0. Analogously $\sum{lk(\triangle_{e}, AE_2)}=0$.

$\bullet$ The fourth equality is an analogue of Remark \ref{equality} which is proved analogously.

$\bullet$ The last equality is Lemma 3$'$ because by the choice of $\alpha$ the map $f:K_5\to\alpha$ is a general position map.

Hence for some edge $e$ of the graph $K_5$ we have $lk(\triangle_{e}, Ae)=1$. Then by Lemma 1$'$ the cycles $\triangle_{e}$ and $Ae$ are linked. {\it QED}

\begin{Pro}

[BE01, \S1.5] Any two closed broken lines in the plane that are in general position intersect each other at an even number of points.
\end{Pro}
This remark is a generalization of Remark 4. We do not prove it here.

\smallskip

{\it Proof of Lemma 3} $'$. By Example 2$'$ it suffices to prove that $v(f)$ does not depend on the map $f:K_5\to\R^2$.
For this it suffices to prove that if general position maps $f,f':K_5\to\R^2$ are equal
on some subgraph $K_4$ then $v(f)=v(f')$.
Assume that these maps are equal on a subgraph $K_4$. Denote by $V$ the vertex of $K_5$ that is not contained in $K_4$.
For any vertex $A$ of $K_4$ denote by $\triangle_A$ the cycle in $K_4$ on three vertices other than $A$, by $VA$ the edge of $K_5$ joining vertices $V$ and $A$.
For any two maps $g,h:K_5\to \R^2$ that are equal on $K_4$ denote $M_{g,h}:=\{g(VA), h(VA), g(e)|\ A\in V(K_4),\ e\in E(K_4)\}$.

{\it Proof that $v(f)=v(f')$ when $f(V)\not=f'(V)$ and $M_{f,f'}$ is a general position set.}
There exists broken line $L$ with the endpoints $f(V)$ and $f'(V)$ such that the set $\{L\}\cup M_{f,f'}$ is a general position set.
The proof of this case is obtained from the part `proof that $v(f) = v(f')$ when $f\cup \{V'\}$ is a general position set' in the proof of Lemma \ref{vankam} by the following changes: replace `segment $VV'$' by `broken line $L$' and use Remark 4$'$ to prove that general position closed broken lines $f(VA)\cup f'(VA)\cup L$ and $f(\triangle_A)$ intersect at an even number of points instead of using Remark \ref{0-even} in the proof that $|VV'A\cap\triangle_A|$ is even.


{\it Proof that $v(f) = v(f')$ in general.}
There exists a general position map $f'':K_5\to\R^2$ such that

$\bullet$ maps $f,f',f'':K_5\to\R^2$ are equal on $K_4$;

$\bullet$ $f''(V)\not=f(V)$, $f''(V)\not=f'(V)$;

$\bullet$ $M_{f,f''}$ and $M_{f',f''}$ are a general position set.

The previous case implies that $v(f)=v(f'')=(f')$. {\it QED}
\comment
So it suffices to prove the item (b) for the case when the map $f':K_5\to \pi$ satisfies the above conditions for the map $f'':K_5\to\pi$. Assume that the map $f':K_5\to\pi$ satisfy these conditions.
There exists broken line $L$ with endpoints $f(v)$ and $f'(v)$ that is in general position with $f(e)$ and $f'(e)$ for any edge $e$ of the graph $K_5$.
For any vertex $i$ of the graph $K_4$ closed broken lines $f(vi)\cup f'(vi)\cup L$ and $f(\triangle_i)$ are in general position. Then Proposition implies that these broken lines intersect at an even number of points. The proof of the fact that $v(f)=v(f')$
is analogous to the proof of Lemma \ref{vankam}. Only the following changes should be made: replace `segment $f(v)f'(v)$' with `broken line $L$' and use Proposition instead of the fact that any two triangles whose vertices are in general position in the plane intersect at an even number of points. {\it QED.}

\endcomment
\bigskip

{\bf Proof of the Sachs Theorem.}

In order to prove the Sachs Theorem we will use the following known lemma.

\begin{lemma3''}\label{vankam3} [Sk]
For any general position map $f:K_{3,3}\to\R^2$ we have $v(f)=1$.


\end{lemma3''}

The proof of this lemma is analogous to the proofs of Lemma \ref{vankam} and Lemma 3$'$.

\smallskip
\textit{Proof of the Sachs Theorem.}
Assume that the graph $K_{4,4}$ is piecewise-linearly embedded in space.
By Lemma \ref{gen} there exists a plane $\alpha$ in general position to $K_{4,4}$.
Denote by $A,B$ some two vertices of $K_{4,4}$ from different parts. Denote $K_{3,3}:=K_{4,4}-A-B$. Let $f:K_{3,3} \to \alpha$ be the orthogonal projection.
For any edge $e$ of $K_{3,3}$ denote by 

$\bullet$ $E_1,E_2$ the vertices joined by $e$ so that $E_1,A$ are in one part, $E_2,B$ are in the other one;

$\bullet$ by $ABe$ the cycle in $K_{4,4}$ on the vertices $A,B,E_1,E_2$; 

$\bullet$ by $\triangle_{e}$ the cycle in $K_{3,3}$ on four vertices other than $E_1$ and $E_2$.

From now on, in any sum, if the limits of the summation are not written, the sum is over edges $e$ of the graph $K_{3,3}$. We have

$$\sum{lk(\triangle_{e}, ABe)}\overset{1}=\sum{lk(\triangle_{e}, AB)}+\sum{lk(\triangle_{e}, AE_2)}+\sum{lk(\triangle_{e}, BE_1)}+\sum{lk(\triangle_{e}, e)}\overset{2}=$$
$$=\sum{lk(\triangle_{e}, e)}\overset{3}=\sum\limits_{(e',e)\in \widetilde{E}(K_{3,3})}{lk(e',e)}\overset{4}= v(f)\overset{5}= 1 \mod 2.$$

$\bullet$ First and third equalities are clear.

$\bullet$ Let us prove the second equality. For any edge $e$ of $K_{3,3}$ we have $lk(\triangle_{e}, AB)=\sum\limits_{e'\in E(\triangle_{e})}{lk(e', AB)}$. Any edge $e'$ of $K_{3,3}$ is contained in exactly four cycles of length 4 in $K_{3,3}$. Then the number $lk(e', AB)$ `appears four times' in the sum
$\sum{lk(\triangle_{e}, AB)}$. Hence this sum is equal to 0.
For any edge $e$ of $K_{3,3}$ we have $lk(\triangle_{e}, AE_2)=\sum\limits_{e'\in E(\triangle_e)}{lk(e', AE_2)}$. Any edge $e'$ of the graph $K_{3,3}-E_2$ is contained in exactly two cycles of length 4 in $K_{3,3}-E_2$. Then the number $lk(e', AE_2)$ `appears twice' in the sum $\sum{lk(\triangle_{e}, AE_2)}$. Hence this sum is equal to 0. Analogously $\sum{lk(\triangle_{e}, BE_1)}=0$.

$\bullet$ The fourth equality is an analogue of Remark \ref{equality} which is proved analogously.

$\bullet$ The last equality is Lemma 3$''$ because by the choice of $\alpha$ the map $f:K_{3,3}\to\alpha$ is a general position map.

Hence for some edge $e$ of the graph $K_{3,3}$ we have $lk(\triangle_{e}, ABe)=1$. Then by Lemma 1$'$ the cycles $\triangle_{e}$ and $ABe$ are linked. {\it QED}

\bigskip

The author is grateful to Arkady Skopenkov and Mikhail Skopenkov for productive discussions. The author is also grateful to anonymous referees of Moscow Mathematical Conference for High-School Students for their comments.


\begin{thebibliography}{RSS95}
\bibitem[Sk]{Sk}
\emph{A. Skopenkov.} Algorithms for recognition of the realizability of hypergraphs, in Russian.
See \texttt{http://www.mccme.ru/circles/oim/algor.pdf}.
\bibitem[BE01]{BE01}
\emph{V.G. Boltyansky, V.A. Efremovich.} Intuitive Combinatorial Topology. Springer; 2001 edition. Russion translation: \texttt{http://ilib.mccme.ru/djvu/geometry/boltiansky-nagl-topo.htm}.
\bibitem[PS05]{PS}
\emph{V. Prasolov, M. Skopenkov.} Ramsey theory of knots and links, Mat. Prosv. 3rd ser. 9 (2005), 108-115 (in Russian).
\bibitem[P06]{P06}
\emph{V. Prasolov.} Elements of combinatorial and differential topology, American Mathematical Soc., 2006. Russian translation: \texttt{ftp://ftp.mccme.ru/users/prasolov/topology/topol2.pdf}.
\bibitem[SS14]{SS14}
\emph{A. Skopenkov and M. Skopenkov.} Some short proofs of the nonrealizability of hypergraphs. See \texttt{http://arxiv.org/pdf/1402.0658.pdf}.
\end{thebibliography}
\end{document}